# The Bias and Efficiency of Incomplete-Data Estimators in Small Univariate Normal Samples

Paul T. von Hippel

LBJ School of Public Affairs, University of Texas
2515 Red River, Box Y, Austin, TX 78712

(512) 537-8112

paulvonhippel.utaustin@gmail.com

## *Abstract*

Widely used methods for analyzing missing data can be biased in small samples. To understand these biases, we evaluate in detail the situation where a small univariate normal sample, with values missing at random, is analyzed using either observed-data maximum likelihood (ML) or multiple imputation (MI). We evaluate two types of MI: the usual Bayesian approach, which we call posterior draw (PD) imputation, and a little-used alternative, which we call ML imputation, in which values are imputed conditionally on an ML estimate. We find that observed-data ML is more efficient and has lower mean squared error than either type of MI. Between the two types of MI, ML imputation is more efficient than PD imputation, and ML imputation also has less potential for bias in small samples. The bias and efficiency of PD imputation can be improved by a change of prior.

4600 words, 2 Appendices, 3 Tables, 3 Figures.

# 1  INTRODUCTION

After many years of treating missing data in an *ad hoc* fashion, analysts are now turning to two more principled methods. One method is observed-data maximum likelihood (ML); the other method is multiple imputation (MI) (Little and Rubin 2002; Allison 2002). There are several forms of MI. By far the most popular form—which we call *posterior draw (PD) imputation*—imputes values by drawing from the Bayesian posterior predictive distribution (Rubin 1987). A somewhat neglected form of MI—which we call we call *ML imputation*—draws imputations conditionally on an ML estimate, or on another estimate with the same asymptotic standard error as an ML estimate (Wang and Robins 1998).

The large-sample properties of these methods are excellent. If the methods employ the same, correctly specified statistical model, and if values are missing at random, then PD imputation, ML imputation, and observed-data ML are all consistent. In fact, as the number of observations and the number of imputations increase toward infinity, PD imputation, ML imputation and observed-data ML all converge to the same asymptotic point estimate with the same asymptotic standard error (Wang and Robins 1998). However, observed-data ML converges more quickly than ML imputation, and ML imputation converges more quickly than PD imputation (Wang and Robins 1998).

The small-sample properties of missing-data estimators are not as well understood. Simulations have found that PD imputation can suffer from bias and inefficiency in small samples of normal data, even when the imputation model is correctly specified (e.g., Demirtas, Freels, and Yucel 2008; Hoogendoorn and Allison 2009). Observed-data ML can also have small-sample biases, although in simulations the biases of ML have so far been milder than the biases of PD imputation (Yuan, Wallentin, and Bentler 2012). Although simulations have demonstrated that PD imputation and observed-data ML can have small-sample biases, the underlying reasons for those biases are not clear. As for ML imputation, as far as we know its small sample properties not been evaluated.

In this paper, we examine in detail the simple situation where a small sample of univariate normal data is missing values at random. This setting has limited practical interest to data analysts, but the simplicity of univariate settings has made them a recurring theoretical proving ground in the missing data literature (e.g., Rubin and Schenker 1986; Wang and Robins 1998; Horton, Lipsitz, and Parzen 2003; He and Raghunathan 2006). Estimators that perform poorly in such a simple setting seem unlikely to improve if the situation gets more complicated. And in a simple setting it is possible to derive not just the asymptotic properties but the exact small-sample distribution of competing estimators.

In our evaluation, we find that PD imputation, ML imputation, and observed-data ML estimation all have potential for bias and inefficiency in small samples. The biases are limited to estimation of the variance $\sigma^2$ and standard deviation $\sigma$; estimates of the mean $\mu$ are unbiased but can sometimes be inefficient. Of the three methods, PD imputation has the most potential for bias. The bias of PD imputation originates in the extra variation that is added when a posterior draw is taken, especially if the Bayesian prior is diffuse. A better choice of prior can improve the bias and efficiency of PD imputation. The small-sample properties of ML imputation and observed-data ML can also be improved by small modifications.

Both before and after these adjustments, we find that ML has smaller mean squared error than ML imputation, which in turn has smaller mean squared error than PD imputation. This result has been proved in general for large samples (Wang and Robins 1998). Ours is the first study to find that, at least in the univariate normal setting, the result holds for small samples as well.

## 2   POINT ESTIMATES

We will derive the exact distribution of each estimator and then calculate each estimator's expectation, bias, standard error (SE) and root mean square error (RMSE). The calculations are straightforward though sometimes tedious. The tedium was relieved by Mathematica software, version 8.

Suppose that we have a simple random sample of $n$ values from an infinite population of a normal variable $Y$ with mean $\mu$ and variance $\sigma^2$. Of these $n$ values, $n_{obs}$ values are observed, and $n_{mis}=n-n_{obs}$ values are missing. We assume that values are *missing at random* (MAR), which in the univariate setting also means that they are *missing completely at random* (MCAR) (Rubin 1976; Heitjan and Basu 1996). That is, the $n_{obs}$ observed $Y$ values are selected at random from the $n$ sampled cases, so the observed values are a random sample from the population.

In the next section, we obtain the exact bias and efficiency of two different kinds of observed-data estimators: ML-like estimators and PD estimators. After that we show what happens when happens when these observed-data estimates are used to impute missing values.

### 2.1   Observed-data estimators

The distribution of the observed values can be summarized as follows:

$$Y_{obs,i} = \mu + \sigma Z_{obs,i}, \text{ where } Z_{obs,i} \sim N(0,1), i = 1, \ldots, n_{obs} \tag{1}$$

Since the observed values are a random sample from the population, we can obtain consistent estimates from the observed values alone. These are the *observed-data estimates* $\hat{\mu}_{obs}, \hat{\sigma}_{obs}^2, \hat{\sigma}_{obs}$. Several observed-data estimators are available.

### 2.1.1 ML and ML-like estimators

The simplest observed-data estimators are the mean, variance, and standard deviation of the $n_{obs}$ observed values. The sampling distribution of these estimators is familiar, but a brief review will serve to introduce notation. The mean $\hat{\mu}_{obs,M}$ of the observed values has a normal sampling distribution.

$$\hat{\mu}_{obs,M} = \bar{Y}_{obs} = \mu + \sigma\bar{Z}_{obs}$$

$$\text{where } \bar{Y}_{obs} = \frac{1}{n_{obs}}\sum_{i=1}^{n_{obs}} Y_i \text{ and } \bar{Z}_{obs} = \frac{1}{n_{obs}}\sum_{i=1}^{n_{obs}} Z_i \sim N\left(0, \frac{1}{n_{obs}}\right) \tag{2}$$

The variance $\hat{\sigma}_{obs,M}^2$ and standard deviation $\hat{\sigma}_{obs,M}$ of the observed values rely on the centered sum of squares (CSS), which has a scaled chi-square sampling distribution.

$$\hat{\sigma}_{obs,M}^2 = \frac{\text{CSS}}{\nu_{obs} + c_M}, \text{where } \nu_{obs} = n_{obs} - 1$$

$$\text{and } CSS = \sum_{i=1}^{n_{obs}} (Y_i - \bar{Y}_{obs})^2 = \sigma^2 U_{obs} \tag{3}$$

$$\text{and } U_{obs} = \sum_{i=1}^{n_{obs}} (Z_i - \bar{Z}_{obs})^2 \sim \chi_{\nu_{obs}}^2$$

The properties of these estimators depend on the constant $c_M$ which appears in the denominator of $\hat{\sigma}_{obs,M}^2$. If we set the constant $c_M$ to 1, then $\hat{\mu}_{obs,M}$, $\hat{\sigma}_{obs,M}^2$, and $\hat{\sigma}_{obs,M}$ are ML estimators, which are consistent and efficient in large samples, but biased in small samples for $\sigma^2$ and $\sigma$. Choosing a different value of $c_M$ yields *ML-like* estimators with the same large-sample distribution but different small-sample properties. For example, choosing $c_M = 0$ makes $\hat{\sigma}_{obs,M}^2$ the minimum variance unbiased (MVU) estimator for $\sigma^2$, which is unbiased but has a larger SE and RMSE than the ML estimator. And choosing $c_M = 2$ makes $\hat{\sigma}_{obs,M}^2$ the minimum mean square error (MMSE) estimator of $\sigma^2$, which is more biased than the ML estimator but has a smaller SE and the smallest possible RMSE (Theil and Schweitzer 1961).

Choices of $c_M$ that are optimal for estimating $\sigma^2$ are not optimal for estimating $\sigma$. For example, the value $c_M = 0$, which yields an unbiased estimate of $\sigma^2$, yields a negatively biased estimate $\sigma$. If we want the estimator $\hat{\sigma}_{obs,M}$ to be approximately unbiased, we should choose $c_M \approx -1.5$. If we want $\hat{\sigma}_{obs,M}$ to have minimal RMSE, we should choose $c_M \approx -.5$.

Table 1a justifies these assertions by giving formulas for bias and SE, which were obtained by applying expectations to the distributions of $\hat{\mu}_{obs,M}$, $\hat{\sigma}_{obs,M}^2$, and $\hat{\sigma}_{obs,M}$.

It is worth remembering that the small-sample bias and efficiency of $\hat{\sigma}_{obs,M}^2$ and $\hat{\sigma}_{obs,M}$ depend on the constant $c_M$. We will shortly find that the bias and efficiency of the PD estimators also depend on a single constant.

### 2.1.2 PD estimators

*Posterior draw (PD) estimators* are drawn at random from the posterior distribution of the parameters given the observed values. In the missing data literature it is popular to assume that the prior is in some sense noninformative. For the mean $\mu$, the usual noninformative prior is uniform. For the variance $\sigma^2$ some commonly used noninformative priors share the form

$$f_{prior}(\sigma^2) \propto \sigma^{-\nu_{prior}-2} \tag{4}$$

which is a degenerate form of the scaled inverse chi-square distribution with $\nu_{prior}$ degrees of freedom (Schafer 1997; Kim 2004). The bias and efficiency of the PD estimates depends on the value chosen for $\nu_{prior}$. Nearly all the imputation literature chooses $-2 \leq \nu_{prior} \leq 0$ (Rubin and Schenker 1986; Demirtas et al. 2008; Schafer 1997; StataCorp 2009), but, as we will see, better choices are available.

The PD variance estimate $\hat{\sigma}_{obs,PD}^2$ is obtained by dividing the CSS by a random chi-square variable (Rubin and Schenker 1986; Kim 2004):

$$\hat{\sigma}_{obs,PD}^2 = \frac{CSS}{U_{PD}}$$
$$\text{where } U_{PD} \sim \chi_{\nu_{PD}}^2 \text{ and } \nu_{PD} = \nu_{prior} + \nu_{obs} \tag{5}$$

And the PD mean estimator $\hat{\mu}_{obs,PD}$ is drawn from a normal distribution whose standard deviation depends on $\hat{\sigma}_{obs,PD}$ and on the sample size $n_{obs}$:

$$\hat{\mu}_{obs,PD} = \hat{\mu}_{obs,M} + \hat{\sigma}_{obs,PD} Z_{PD}, \text{where } Z_{PD} \sim N\left(0, \frac{1}{n_{obs}}\right) \tag{6}$$

By plugging the distribution of $\bar{Z}_{obs}$ and CSS into the definitions of the PD estimators, we can calculate the distribution of the PD estimators. $\hat{\sigma}_{obs,PD}^2$ is the scaled ratio of two chi-square variables, and therefore follows a scaled $F$ distribution.

$$\hat{\sigma}^2_{obs,PD} = \frac{\sigma^2 U_{obs}}{U_{PD}} = \sigma^2 \frac{\nu_{obs}}{\nu_{PD}} F_{PD}$$

$$\text{where } F_{PD} = \frac{U_{obs}/\nu_{obs}}{U_{PD}/\nu_{PD}} \sim F_{\nu_{obs}, \nu_{PD}} \tag{7}$$

Therefore $\hat{\sigma}_{obs,PD}$ is proportional to the square root of an $F$ variable, and $\hat{\mu}_{obs,PD}$ is a function of three independent variables—two normal and one $F$.

$$\hat{\mu}_{obs,PD} = \mu + \sigma \left( \bar{Z}_{obs} + \sqrt{\frac{\nu_{obs}}{\nu_{PD}}} \sqrt{F_{PD}} Z_{PD} \right) \tag{8}$$

Now that we have the distributions of the PD estimators, we can derive formulas for their bias and SE by taking expectations and variances. Table 1a gives the resulting formulas. Table 1b compares the bias, SE, and RMSE of the ML-like estimators with $c_M = 0,1,2$ (abbreviated M0, M1, M2) and the PD estimators with $\nu_{prior} = -2,0,\dots,7$ (abbreviated PD–2, PD0, …, PD7) in samples of $n_{obs} = 5,20,$ or $100$ observed values from a normal population with $\mu = \sigma = 1$. Figure 1 plots the RMSE of selected estimators as a function of the observed sample size $n_{obs}$.

#### ←Figure 1 near here→

We remarked earlier that nearly all the imputation literature chooses a prior with $-2 \leq \nu_{prior} \leq 0$ degrees of freedom (Rubin and Schenker 1986; Demirtas et al. 2008; Schafer 1997; StataCorp 2009) but, as Table 1 shows, these choices yield estimates that have positive bias for $\sigma^2$ and large or undefined standard errors for $\mu, \sigma,$ and $\sigma^2$.

There are several better choices for $\nu_{prior}$. If we want to optimize the properties of the variance estimator $\hat{\sigma}^2_{obs,PD}$, $\nu_{prior} = 2$ makes $\hat{\sigma}^2_{obs,PD}$ unbiased (cf. Kim 2004), and $\nu_{prior} \approx 7$ yields an estimator $\hat{\sigma}^2_{obs,PD}$ that, though negatively biased, has minimal RMSE. Nearly as small a RMSE can be achieved with less bias by choosing $\nu_{prior} = 6$.

If instead we want to optimize the properties of the standard deviation estimator $\hat{\sigma}_{obs,PD}$, we can make $\hat{\sigma}_{obs,PD}$ unbiased by choosing $\nu_{prior} = 1$, and we can minimize the RMSE of $\hat{\sigma}_{obs,PD}$ by choosing $\nu_{prior} \approx 4$.

Finally, if we want to optimize the properties of the mean estimator $\hat{\mu}_{obs,PD}$, we should choose the largest value of $\nu_{prior}$ that we can. Increasing $\nu_{prior}$ reduces the SE of $\hat{\mu}_{obs,PD}$, and regardless of $\nu_{prior}$ the estimator $\hat{\mu}_{obs,PD}$, will remain unbiased. However, if $\nu_{prior}$ gets too large the bias of $\hat{\sigma}^2_{obs,PD}$ and $\hat{\sigma}_{obs,PD}$ will become unacceptable.

On balance, a reasonable case can be made for any prior in the range $1 \leq \nu_{prior} \leq 7$. While the range of good options is fairly broad, it does not include the most commonly used values; that is, the good options do not include $-2 \leq \nu_{prior} \leq 0$.

In large samples (e.g., $n_{obs} = 100$), the bias of the PD estimators becomes negligible, but the PD estimators remain much less efficient than the ML-like estimators. Specifically, the asymptotic standard errors of the PD estimators are $\sqrt{2}$ times the asymptotic standard errors of the ML estimators. In other words, using a PD estimator is asymptotically equivalent to discarding half the observed values and applying an ML-like estimator to the values that remain. In a sense the inefficiency of the PD estimators is obvious from the way they are constructed; the PD estimators start with ML-like estimators and then add random variables that double the variance (equations (7) and (8)).

## 2.2   Imputation-based estimates

Having obtained observed-data estimates $\hat{\mu}_{obs}, \hat{\sigma}_{obs}$, we can plug them into an equation that imputes the missing values:

$$Y_{imp,i} = \hat{\mu}_{obs} + \hat{\sigma}_{obs} Z_{imp,i}, \text{where } Z_{imp,i} \sim N(0,1), i = n_{obs} + 1, \dots, n \qquad (9)$$

If the imputation model uses the ML-like observed-data estimators $\hat{\mu}_{obs,M}, \hat{\sigma}_{obs,M}$, we call this process ML imputation. If the imputation model uses the PD observed-data estimators $\hat{\mu}_{obs,PD}, \hat{\sigma}_{obs,PD}$, we call it PD imputation. If we impute the data once, we call the process single imputation (SI); if we impute the data several times we call the process multiple imputation (MI). We will begin by deriving the properties of SI estimates, and then move on to MI estimates.

### 2.2.1   Single imputation (SI) estimators

The SI sample $Y_{SI} = (Y_{obs}, Y_{imp})$ is a mixture of two normal distributions. The mixture contains $n_{obs}$ observed values $Y_{obs}$ drawn at random from a normal population with mean $\mu$ and variance $\sigma^2$, mixed with $n_{mis}$ imputed values $Y_{imp}$ drawn at random from a slightly different normal population with mean $\hat{\mu}_{obs}$ and variance $\hat{\sigma}_{obs}^2$.

After imputation, we re-estimate $\mu$, $\sigma^2$, and $\sigma$ using the mean, variance, and standard deviation of the SI sample. This results in the following SI estimators:

$$\hat{\mu}_{SI} = \bar{Y}_{SI} = \frac{1}{n} \sum_{i=1}^{n} Y_{SI,i}$$

$$\hat{\sigma}_{SI}^2 = s_{SI}^2 = \frac{1}{n-1} \sum_{i=1}^{n} (Y_{SI,i} - \bar{Y}_{SI})^2 \qquad (10)$$

$$\hat{\sigma}_{SI} = \sqrt{\hat{\sigma}_{SI}^2}$$

To obtain the distribution of the SI mean $\hat{\mu}_{SI}$, we start by breaking it into two components—the mean of the observed values and the mean of the imputed values:

$$\hat{\mu}_{SI} = \frac{n_{obs} \bar{Y}_{obs} + n_{mis} \bar{Y}_{imp}}{n} \qquad (11)$$

Likewise, to calculate the distribution of the SI variance $\hat{\sigma}_{SI}^2$, we break $\hat{\sigma}_{SI}^2$ into three components—the variance $s_{obs}^2$ within the observed values, the variance $s_{imp}^2$ within the imputed values, and the variance $s_{btw}^2$ between the observed values and the imputed values:

$$\hat{\sigma}_{SI}^2 = \frac{1}{n-1} \left( (n_{obs}-1)s_{obs}^2 + (n_{mis}-1)s_{imp}^2 + s_{btw}^2 \right) \qquad (12)$$

We obtain the distributions of $\hat{\mu}_{SI}$, $\hat{\sigma}_{SI}^2$, and $\hat{\sigma}_{SI}$ by plugging in the distribution of the observed-data summaries $\bar{Y}_{obs}$ and $s_{obs}^2$, along with the distributions of $\bar{Y}_{imp}$, $s_{imp}^2$, and $s_{btw}^2$, under ML imputation and under PD imputation. The distributions of $\hat{\mu}_{SI}$, $\hat{\sigma}_{SI}^2$, and $\hat{\sigma}_{SI}$ are derived in Appendix SI.

Once we have the distributions of the SI estimators, we derive their bias and SE by taking expectations and variances. Table 2a gives the resulting formulas. Table 2b illustrates the results by displaying the expectation, bias, SE, and RMSE of the SI estimators for a normal variable with $\mu = \sigma = 1$ and sample sizes of $n_{obs} = n_{mis} = 5, 20,$ and $100$. The numeric values in Table 2b were calculated directly from the formulas in Table 2a, and then verified by simulation.

<center>←Table 2 <strong>near here</strong>→</center>

In some ways the results for the SI estimators are similar to those for the observed-data estimators. Under SI, the PD estimators have larger RMSEs than the ML-like estimators. And the most popular PD estimators (PD0 and PD–2) are the worst, with large biases and low efficiency when estimating $\sigma^2$ in small samples.

But the differences among the SI estimators are relatively small compared to the differences among the observed-data estimators. This is because the SI estimators combine information across the imputed values $Y_{imp}$ and the observed values $Y_{obs}$. Combining $Y_{imp}$ and $Y_{obs}$ increases the RMSEs of the ML-like estimators by adding random variation from $Y_{imp}$, but it reduces the RMSEs of the PD estimators by smoothing them toward the efficient observed-data estimators $\bar{Y}_{obs}$, $s^2_{obs}$, and $s_{obs}$. By adding variation to the efficient ML-like estimators, while reducing variation in the less-efficient PD estimators, the SI process brings all of the estimators closer together.

### 2.2.2   Multiple imputation (MI) estimators

We now progress from SI to MI. MI is an iterative process. In iteration $d = 1, \ldots, D$, we carry out the following steps:

1. Calculate observed-data point estimates $\hat{\mu}_{obs,d}$ and $\hat{\sigma}^2_{obs,d}$.
2. Impute random values conditionally on the observed-data estimates, yielding SI data.
3. Analyze the SI data to obtain SI point estimates $\hat{\mu}_{SI,d}$, $\hat{\sigma}^2_{SI,d}$, $\hat{\sigma}_{SI,d}$.

Under PD imputation, different observed-data PD estimates $\hat{\mu}_{obs,PD,d}$ and $\hat{\sigma}^2_{obs,PD,d}$ are drawn in every iteration, so all three steps of the algorithm must be iterated. But under ML imputation, the observed-data ML-like estimates $\hat{\mu}_{obs,M}$, $\hat{\sigma}^2_{obs,M}$ are the same in every iteration, so we can run step 1 once and then iterate steps 2 and 3.

After repeating these steps $D$ times, we average the $D$ SI point estimates to obtain MI point estimates. For our purposes it will be helpful to imagine *infinite imputation* ($\infty$I) estimators, which limit the MI estimators as the number of imputations $D$ increases.

$$\hat{\mu}_{MI} = \frac{1}{D} \sum_{m=1}^{D} \hat{\mu}_{SI,d} \xrightarrow[D \to \infty]{} \hat{\mu}_{\infty I}$$

$$\hat{\sigma}^2_{MI} = \frac{1}{D} \sum_{d=1}^{D} \hat{\sigma}^2_{SI,d} \xrightarrow[D \to \infty]{} \hat{\sigma}^2_{\infty I} \qquad (13)$$

$$\hat{\sigma}_{MI} = \frac{1}{D} \sum_{d=1}^{D} \hat{\sigma}_{SI,d} \xrightarrow[D \to \infty]{} \hat{\sigma}_{\infty I}$$

Averaging across several imputations does not change the bias of imputation-based estimators. The expectation, and therefore the bias, is the same whether the number of imputations is one (SI), several (MI), or infinite ($\infty$I):

$$E(\hat{\mu}_{SI}) = E(\hat{\mu}_{MI}) = E(\hat{\mu}_{\infty I})$$
$$E(\hat{\sigma}^2_{SI}) = E(\hat{\sigma}^2_{MI}) = E(\hat{\sigma}^2_{\infty I}) \qquad (14)$$
$$E(\hat{\sigma}_{SI}) = E(\hat{\sigma}_{MI}) = E(\hat{\sigma}_{\infty I})$$

Although averaging across imputations does not affect bias, it does increase efficiency. As the number of imputations grows, the standard error of the MI estimator approaches the standard error of the $\infty$I estimator. This observation is helpful since it permits us to calculate the efficiency of the MI estimators. Note that, with a finite number of imputations $D$, the variance of the MI estimator is a weighted average of the variance of the single imputation estimator and the variance of the infinite imputation estimator. That is, for an MI estimate of $\mu$ the variance is

$$
\begin{aligned}
V(\hat{\mu}_{MI}) &= V(\hat{\mu}_{\infty I}) + V(\hat{\mu}_{MI}|\hat{\mu}_{\infty I}) \\
&= V(\hat{\mu}_{\infty I}) + \frac{1}{D}V(\hat{\mu}_{SI}|\hat{\mu}_{\infty I}) \\
&= V(\hat{\mu}_{\infty I}) + \frac{1}{D}\big(V(\hat{\mu}_{SI}) - V(\hat{\mu}_{\infty I})\big) \\
&= \Big(1 - \frac{1}{D}\Big)V(\hat{\mu}_{\infty I}) + \frac{1}{D}V(\hat{\mu}_{SI})
\end{aligned}
\tag{15}
$$

And for MI estimates of $\sigma$ and $\sigma^2$ the variance is

$$
\begin{aligned}
V(\hat{\sigma}_{MI}) &= \Big(1 - \frac{1}{D}\Big)V(\hat{\sigma}_{\infty I}) + \frac{1}{D}V(\hat{\sigma}_{SI}) \\
V(\hat{\sigma}_{MI}^2) &= \Big(1 - \frac{1}{D}\Big)V(\hat{\sigma}_{\infty I}^2) + \frac{1}{D}V(\hat{\sigma}_{SI}^2)
\end{aligned}
\tag{16}
$$

Now we use these formulas to calculate the variance of an MI estimator. We begin by deriving the variance of the $\infty$I estimators, which we obtain by taking the variance of the expectations of the SI estimators, conditionally on the observed values:

$$
\begin{aligned}
V(\hat{\mu}_{\infty I}) &= V(E(\hat{\mu}_{SI}|Y_{obs})) \\
V(\hat{\sigma}_{\infty I}) &= V(E(\hat{\sigma}_{SI}|Y_{obs})) \\
V(\hat{\sigma}_{\infty I}^2) &= V(E(\hat{\sigma}_{SI}^2|Y_{obs}))
\end{aligned}
\tag{17}
$$

To calculate these conditional variances, Appendix $\infty$I derives the conditional *distributions* of the $\infty$I estimators by taking conditional expectations—

$$
\begin{aligned}
f(\hat{\mu}_{\infty I}) &= E(\hat{\mu}_{SI}|Y_{obs}) \\
f(\hat{\sigma}_{\infty I}) &= E(\hat{\sigma}_{SI}|Y_{obs}) \\
f(\hat{\sigma}_{\infty I}^2) &= E(\hat{\sigma}_{SI}^2|Y_{obs})
\end{aligned}
\tag{18}
$$

Given the distributions of the $\infty$I estimators, it is straightforward to calculate their variances. These variances, or actually the corresponding SEs, are given in Table 3a.

Now we obtain the variances of the MI estimators by using equations (15) and (16) to combine the variances of the $\infty$I estimators (in Table 3a) with the variances of the SI estimators (Table 2a). Table 3a gives the resulting formulas for the SEs of the MI estimators. The biases of the MI estimators are not given because they are the same as the biases of the SI estimators, which were given in Table 2b.

Table 3b illustrates the results by calculating the expectation, bias, SE, and RMSE for MI estimators with $D = 5$ imputations in samples with $n_{obs} = n_{mis} = 5, 20$, and $100$ observations from a population with $\mu = \sigma = 1$. The estimators are evaluated under ML imputation with $c_M = 0,1,2$ (abbreviated M0, M1, M2) and under PD imputation with $v_{prior} = -2,0,...,7$ (abbreviated PD–2, PD0, ..., PD7). The values in Table 3b were calculated directly from the formulas in Table 3a, and then verified by simulation.

Figure 2 and Figure 3 plot the RMSE of each MI estimator as a function of the number of imputations $D$ and the number of observed values $n_{obs}$, under the assumption that half of values are missing (i.e., $n_{mis} = n_{obs}$). In Figure 2, $n_{obs}$ increases while $D$ is held constant at 5. In Figure 3, $D$ increases while $n_{obs}$ is held constant at 20.

Under PD imputation, the popular PD0 and PD–2 estimators give the worst estimates with large or undefined RMSEs and large or undefined biases in small samples. The PD2, PD4, and PD6 estimators give much better estimates, with RMSEs that are often comparable to the RMSEs under ML imputation.

In very small samples, the RMSEs of different estimators can be quite different, but the estimators converge as $D$ increases and, especially, as $n_{obs}$ increases. The rate of convergence depends on the parameter being estimated. MI estimators of the mean $\mu$ converge faster than estimators of the standard deviation $\sigma$, which converge faster than estimators of the variance $\sigma^2$.

The differences among the MI estimators are smaller than the differences among the corresponding SI estimators. You can verify this by comparing the SI estimators in Table 2b to the MI estimators in Table 3b, or by noticing, in Figure 3, how the RMSEs of different MI estimators grow more similar as $D$ increases beyond 1.

The reason that MI estimators are more similar than SI estimators is that averaging across multiple imputations reduces random variation. The more random variation an SI estimator has, the more it improves under MI. So the more variable PD estimators, especially the PD–2 and PD0 estimators, benefit more from MI than the less variable ML-like estimators.

To return to an earlier analogy: using a single observed-data PD estimator is asymptotically equivalent to using an ML-like estimator on half the sample. But under MI, a fresh PD estimate is drawn in every iteration, and that is like selecting a *different* half-sample in each iteration. With enough iterations $D$, we are effectively averaging many half-samples, and we approach a situation where all the observed values in the full sample have made equal contributions to the estimate. That is why PD imputation is almost as efficient as ML imputation when the number of imputations $D$ is large.

In fact, in large samples ($n_{obs} \to \infty$) with infinite imputations ($D \to \infty$), $\infty$I estimators under PD imputation are *identical* to $\infty$I estimators under ML imputation (Wang and Robins 1998). In small samples, however, different $\infty$I estimators may have different biases and different standard errors. Table 2 and Table 3 verify this by showing that the bias and standard error of $\hat{\sigma}_{\infty I}$ and $\hat{\sigma}_{\infty I}^2$ depend on the constants $c_M$ and $\nu_{prior}$, which can be important in small samples. For example, in a small sample with $n_{obs} = n_{mis} = 20$, the estimator $\hat{\sigma}_{MI,PD-2}^2$ will have a bias of 14% and an SE that is at least 27% larger than the SE of $\hat{\sigma}_{MI,PD7}^2$—no matter how many imputations $D$ are used. However, in a larger sample with $n_{obs} = n_{mis} = 100$ and $D$=5, both $\hat{\sigma}_{MI,PD-2}^2$ and $\hat{\sigma}_{MI,PD7}^2$ have approximately the same 2% bias (albeit in opposite directions) and approximately the same SEs.

## 3    STANDARD ERROR ESTIMATES

The paper so far has focused on the properties of point estimators. We now discuss how to estimate standard errors.

### 3.1    Observed-data ML-like estimators

The standard error of the observed-data ML-like estimator $\hat{\mu}_{obs,M}$ is estimated by the familiar expression

$$\sqrt{\hat{V}(\hat{\mu}_{obs,M})} = \frac{\hat{\sigma}_{obs,M}}{\sqrt{n_{obs}}} \tag{19}$$

This expression applies to the ML estimator (with $c_M = 0$) and to the ML-like MVU and MMSE estimators (with $c_M = 1$ or 2).

Comparison with the true standard error in Table 1 shows that the estimated standard error $\sqrt{V(\hat{\mu}_{obs,M})}$ is simply the true standard error with an estimate $\hat{\sigma}_{obs,M}$ substituted for $\sigma$. It follows that $\sqrt{V(\hat{\mu}_{obs,M})}$ inherits the properties of $\hat{\sigma}_{obs,M}$. For example, if we use the estimator $\hat{\sigma}_{obs,M}$ with $c_M = 0$ or 1, then in small samples $\hat{\sigma}_{obs,M}$ is negatively biased, so $\sqrt{V(\hat{\mu}_{obs,M})}$ is negatively biased as well. On the other hand, if we choose $c_M \approx -1.5$ then $\hat{\sigma}_{obs,M}$ and $\sqrt{V(\hat{\mu}_{obs,M})}$ are approximately unbiased, and if we choose $c_M \approx -.5$ then $\hat{\sigma}_{obs,M}$ and $\sqrt{V(\hat{\mu}_{obs,M})}$ have minimal RMSE.

The usual practice, of course, is to choose a value of $c_M$ with an eye toward the properties of the point estimators $\hat{\sigma}_{obs,M}$ or $\hat{\sigma}^2_{obs,M}$. The bias and efficiency of the standard error estimate are secondary concerns.

## 3.2   ML imputation

The following formula gives a consistent estimator for the standard error of a scalar estimand $\theta$ that is estimated by ML imputation (Wang and Robins 1998):

$$\sqrt{\hat{V}\left(\hat{\theta}_{MI,M}\right)} = \sqrt{\hat{V}\left(\hat{\theta}_{obs,M}\right) + \frac{\gamma}{D}\,\overline{W}_{MI,M}}$$

$$\text{where } \overline{W}_{MI,M} = \frac{1}{D}\sum_{d=1}^{D}\hat{W}_{SI,M} \qquad (20)$$

Here $\sqrt{\hat{V}\left(\hat{\theta}_{obs,M}\right)}$ is the estimated standard error of the observed-data M estimate, and $\sqrt{\hat{W}_{SI,M}}$ is the standard error estimate that would apply if we treated each SI data set as though all the values were observed. $\gamma$ is the *fraction of missing information*, which in the univariate setting is just the fraction of values that are missing.

So the components of $\hat{V}\left(\hat{\mu}_{MI,M}\right)$ are

$$\hat{V}\left(\hat{\mu}_{obs,M}\right) = \frac{\hat{\sigma}^2_{obs,M}}{n_{obs}}$$

$$\gamma = \frac{n_{mis}}{n} \qquad (21)$$

$$\overline{W}_{MI,M} = \frac{\hat{\sigma}^2_{MI,M}}{n}$$

and the estimated standard error is

$$\sqrt{\hat{V}\left(\hat{\mu}_{MI,M}\right)} = \sqrt{\frac{\hat{\sigma}^2_{obs,M}}{n_{obs}} + \frac{n_{mis}}{Dn^2}\,\hat{\sigma}^2_{MI,M}} \qquad (22)$$

Notice that, if $c_M = 0$, $\sqrt{V(\hat{\mu}_{MI,M})}$ is just the true standard error from Table 3 with estimates $\hat{\sigma}^2_{obs,M}$ and $\hat{\sigma}^2_{MI,M}$ substituted for $\sigma^2$. If $0 \le c_M \le 2$, then $\hat{\sigma}_{obs,M}$ and $\hat{\sigma}_{MI,M}$ have negative bias, so $\sqrt{V(\hat{\mu}_{MI,M})}$ has negative bias as well, but the bias shrinks as the observed sample size $n_{obs}$ grows. Notice also that $\sqrt{V(\hat{\mu}_{MI,M})}$ converges to $\sqrt{V(\hat{\mu}_{obs,M})}$ as $D$ increases.

The practical difficulty posed formula (20) is that the MI standard error estimate $\sqrt{\hat{V}(\hat{\theta}_{MI,M})}$ requires the observed-data standard error estimate $\sqrt{\hat{V}(\hat{\theta}_{obs,M})}$. This is no problem if $\sqrt{\hat{V}(\hat{\theta}_{obs,M})}$ was calculated along with the observed-data estimate $\hat{\theta}_{obs,M}$, but it is often the case that $\sqrt{\hat{V}(\hat{\theta}_{obs,M})}$ was not calculated or, if it was calculated, it was not retained.

Wang and Robins (1998) have proposed an alternative formula for the standard error of an MI estimate under ML imputation. Their formula requires more computational effort, though, and it has not been implemented in available software.

### 3.3 PD imputation

Under PD imputation, the standard error is typically estimated using a formula that relies on summary statistics calculated from the imputed data sets alone (Rubin 1987). For a scalar estimand $\theta$, the formula for the standard error of the PD MI point estimator is based on the variances within and between the $D$ imputed data sets:

$$\sqrt{\hat{V}(\hat{\theta}_{MI,PD})} = \sqrt{\overline{W}_{MI,PD} + \frac{D+1}{D}\hat{B}_{MI,PD}}$$

$$\text{where } \overline{W}_{MI,PD} = \frac{1}{D}\sum_{d=1}^{D}\hat{W}_{SI,PD} \tag{23}$$

$$\text{and } \hat{B}_{MI,PD} = \frac{1}{D-1}\sum_{d=1}^{D}\left(\hat{\theta}_{SI,PD,d} - \hat{\theta}_{MI,PD}\right)^2$$

The within variance $\overline{W}_{MI,PD}$ is the mean square of the standard error estimate that would be appropriate if we analyzed each SI data set as though all its values were observed and none were imputed. The between variance $\hat{B}_{MI,PD}$ is an unbiased estimate of

$$E(\hat{B}_{MI,PD}) = V(\hat{\theta}_{SI,PD}|\hat{\theta}_{\infty I,PD}) = V(\hat{\theta}_{SI,PD}) - V(\hat{\theta}_{\infty I,PD}) \tag{24}$$

So if we estimate the mean $\mu$ of univariate normal data,

$$\hat{W} = \frac{\hat{\sigma}_{SI,PD,d}^2}{n} \text{ so that } \overline{W} = \frac{\hat{\sigma}_{MI,PD}^2}{n}$$

$$\hat{B} = \frac{1}{D-1}\sum_{d=1}^{D}\left(\hat{\mu}_{SI,PD,d} - \hat{\mu}_{MI,PD}\right)^2 \tag{25}$$

and the estimated standard error is

$$\sqrt{\hat{V}(\hat{\mu}_{MI,PD})} = \sqrt{\frac{\hat{\sigma}_{MI}^2}{n} + \frac{D+1}{D(D-1)}\sum_{d=1}^{D}(\hat{\mu}_{SI,PD,d} - \hat{\mu}_{MI,PD})^2} \qquad (26)$$

Using results in this paper, we now show that $\hat{V}(\hat{\mu}_{MI,PD})$, like $\hat{\sigma}_{MI,PD}^2$, is biased unless $v_{prior} = 2$.:

$$
\begin{aligned}
Bias\left(\hat{V}(\hat{\mu}_{MI,PD})\right) &= E\left(\hat{V}(\hat{\mu}_{MI,PD})\right) - V(\hat{\mu}_{MI,PD}) \\
&= E\left(\frac{\hat{\sigma}_{MI,PD}^2}{n} + \frac{D+1}{D}\hat{B}\right) - V(\hat{\mu}_{MI,PD}) \\
&= \frac{E(\hat{\sigma}_{MI,PD}^2)}{n} + \frac{D+1}{D}\left(V(\hat{\mu}_{SI,PD}) - V(\hat{\mu}_{\infty I,PD})\right) - V(\hat{\mu}_{MI,PD}) \\
&= -\sigma^2 \frac{n_{\text{mis}}(n + n_{\text{obs}} - 1)(v_{\text{prior}} - 2)}{(n-1)n n_{\text{obs}}(n_{\text{obs}} + v_{\text{prior}} - 3)} \\
&= \frac{n + n_{\text{obs}} - 1}{n n_{\text{obs}}} Bias(\hat{\sigma}_{MI,PD}^2)
\end{aligned}
\qquad (27)
$$

Kim (2004) reaches a similar conclusion by a different route. This is an argument for using $v_{prior} = 2$. However, we would favor a different value of $v_{prior}$ if we wanted to minimize the RMSE of $\hat{V}(\hat{\mu}_{MI,PD})$ or if we wanted to minimize the bias or RMSE of $\sqrt{\hat{V}(\hat{\mu}_{MI,PD})}$.

We prefer to choose $v_{prior}$ with an eye toward the bias and RMSE of the point estimators—$\hat{\mu}_{MI,PD}$, $\hat{\sigma}_{MI,PD}$, or $\hat{\sigma}_{MI,PD}^2$—which as remarked earlier could justify any choice in the range $1 \leq v_{prior} \leq 7$. The properties of the standard error estimator are a secondary concern.

# 4 CONCLUSION

We have considered three different methods for estimating the mean, variance, and standard deviation of incomplete univariate normal data. All three of the methods that we considered—observed-data ML, ML imputation, and PD imputation—have potential for bias and inefficiency in small samples, but all three methods offer ways to reduce the bias or increase the efficiency. Under PD imputation, the small-sample properties of the estimator hinges on the Bayesian prior. The problem of choosing a Bayesian prior can be avoided if we use ML imputation, which has smaller RMSE. However, the standard errors of point estimates are harder to calculate under ML imputation than under PD imputation. A final option is to avoid imputation altogether and derive ML-like estimates from the observed data alone. Observed-data ML estimates have smaller RMSE than any estimate obtained from imputation.

# APPENDIX SI

This appendix derives the distribution of the SI estimators.

To obtain the distribution of the SI mean $\hat{\mu}_{SI}$, we first break it into two components—the mean of the observed values and the mean of the imputed values:

$$\hat{\mu}_{SI} = \frac{n_{obs}\bar{Y}_{obs} + n_{mis}\bar{Y}_{imp}}{n}$$
$$= \frac{1}{n}\left(n_{obs}(\mu + \sigma\bar{Z}_{obs}) + n_{mis}(\hat{\mu}_{obs} + \hat{\sigma}_{obs}\bar{Z}_{imp})\right) \tag{28}$$

where

$$\frac{1}{n_{mis}}\sum_{i=n_{obs}+1}^{n} Z_{imp,i} = \bar{Z}_{imp} \sim N\left(0, \frac{1}{n_{mis}}\right) \tag{29}$$

Now we can plug into equation (28) ML-like observed-data estimators $\hat{\mu}_{obs,M}, \hat{\sigma}_{obs,M}$ to obtain the distribution of $\hat{\mu}_{SI}$ under ML imputation:

$$\hat{\mu}_{SI,M} = \mu + \sigma\left(\bar{Z}_{obs} + \frac{n_{mis}}{n}\sqrt{\frac{U_{obs}}{\nu_{obs} + c_M}}\bar{Z}_{imp}\right) \tag{30}$$

Or we can plug in the PD observed-data estimators $\hat{\mu}_{obs,PD}, \hat{\sigma}_{obs,PD}$ to obtain the distribution of $\hat{\mu}_{SI}$ under PD imputation:

$$\hat{\mu}_{SI,PD} = \mu + \sigma\left(\bar{Z}_{obs} + \frac{n_{mis}}{n}\sqrt{\frac{U_{obs}}{\nu_{PD}}}\, t_{SI}\right) \tag{31}$$

where

$$\frac{(Z_{PD} + \bar{Z}_{imp})}{\sqrt{U_{PD}/\nu_{PD}}} = t_{SI} \sim t\left(0, \frac{1}{n_{obs}} + \frac{1}{n_{mis}}, \nu_{PD}\right)$$

follows a 3-parameter $t$ distribution.

Similarly, to calculate the distribution of the SI variance $\hat{\sigma}_{SI}^2$, we break $\hat{\sigma}_{SI}^2$ into three components—the variance $s_{obs}^2$ within the observed values, the variance $s_{imp}^2$ within the imputed values, and the variance $s_{btw}^2$ between the observed values and the imputed values:

$$\hat{\sigma}_{SI}^2 = \frac{1}{n-1}\left((n_{obs}-1)s_{obs}^2 + (n_{mis}-1)s_{imp}^2 + s_{btw}^2\right) \tag{32}$$

Here $s_{obs}^2$ has a scaled chi-square distribution—

$$s_{obs}^2 = \frac{1}{\nu_{obs}}\sum_{i=1}^{N}(Y_{obs,i}-\bar{Y}_{obs})^2 = \frac{\sigma^2 U_{obs}}{\nu_{obs}}$$

$$\text{where } U_{obs} = \sum_{i=1}^{n_{obs}}(Z_i - \bar{Z}_{obs})^2 \sim \chi_{\nu_{obs}}^2 \text{ and } \nu_{obs} = n_{obs} - 1$$

—while $s_{imp}^2$ has a more complicated distribution that depends in part on the distribution of $\hat{\sigma}_{obs}^2$—

$$s_{imp}^2 = \frac{1}{\nu_{mis}}\sum_{i=1}^{N}(Y_{imp,i}-\bar{Y}_{imp})^2 = \frac{\hat{\sigma}_{obs}^2 U_{imp}}{\nu_{mis}},$$

$$\text{where } U_{imp} = \sum_{i=n_{obs}+1}^{n}\left(Z_{imp,i}-\bar{Z}_{imp}\right)^2 \sim \chi_{\nu_{mis}}^2 \text{ and } \nu_{mis} = n_{mis} - 1$$

—and $s_{btw}^2$ has a distribution that depends in part on the distribution of both $\hat{\mu}_{obs}$ and $\hat{\sigma}_{obs}^2$—

$$s_{btw}^2 = n_{obs}\left(\bar{Y}_{obs}-\bar{Y}_{comp}\right)^2 + n_{mis}\left(\bar{Y}_{imp}-\bar{Y}_{comp}\right)^2$$

$$= \frac{n_{obs}n_{mis}}{n}\left(\hat{\mu}_{obs}-\mu + \hat{\sigma}_{obs}\bar{Z}_{imp} - \sigma\bar{Z}_{obs}\right)^2$$

To calculate the distribution of $\hat{\sigma}_{SI}^2$ under ML imputation, we plug the ML-like observed-data estimators $\hat{\mu}_{obs,M}, \hat{\sigma}_{obs,M}$ into the definition for $\hat{\sigma}_{SI}^2$. The result is

$$\hat{\sigma}_{SI,M}^2 = \sigma^2 \frac{U_{obs}}{n-1}\left(1 + \frac{1}{\nu_{obs}+c_M}\left(U_{imp}+\frac{n_{obs}}{n}n_{mis}\bar{Z}_{imp}^2\right)\right) \tag{33}$$

$$\text{where } n_{mis}\bar{Z}_{imp}^2 \sim \chi_1^2$$

Likewise, to get the distribution of $\hat{\sigma}_{SI}^2$ under PD imputation, we plug the PD observed-data estimators $\hat{\mu}_{obs,PD}, \hat{\sigma}_{obs,PD}$ into the definition for $\hat{\sigma}_{SI}^2$. The result is

$$\hat{\sigma}_{SI,PD}^2 = \sigma^2 \frac{U_{obs}}{n-1}\left(1 + \frac{n_{mis}}{\nu_{PD}}F_{SI}\right) \tag{34}$$

where

$$F_{SI} = \frac{\left(U_{\text{imp}} + \frac{n_{mis}n_{obs}}{n}\left(\bar{Z}_{imp} + Z_{PD}\right)^2\right)/n_{mis}}{U_{PD}/\nu_{PD}} \sim F_{n_{mis},\nu_{PD}} \tag{35}$$
$$\text{since } \frac{n_{mis}n_{obs}}{n}\left(\bar{Z}_{imp} + Z_{PD}\right)^2 \sim \chi_1^2$$

By taking expectations and variances, we can convert these distributions into formulas for bias and SE. The results are given in Table 2.

# APPENDIX ∞I

This appendix derives the distributions of the ∞I estimators. As indicated in equation (18), these distributions are obtained by starting with the SI estimators and taking expectations of the statistics that vary from one SI dataset to another.

Under PD imputation, the distributions of the ∞I estimators are

$$\hat{\mu}_{\infty I,PD} = \mu + \sigma\left(\bar{Z}_{obs} + \frac{n_{mis}}{n}\sqrt{\frac{U_{obs}}{\nu_{PD}}}\ E(t_{SI,m})\right)$$
$$= \mu + \sigma\bar{Z}_{obs}$$

$$\hat{\sigma}^2_{\infty I,PD} = \sigma^2\frac{U_{obs}}{n-1}\left(1 + \frac{n_{mis}}{\nu_{PD}}E(F_{SI,m})\right)$$
$$= \sigma^2\frac{U_{obs}}{n-1}\left(1 + \frac{1}{\nu_{PD}-2}\right)$$

$$\hat{\sigma}_{\infty I,PD} = \sigma\sqrt{\frac{U_{obs}}{n-1}}E\left(\sqrt{1+\frac{n_{mis}}{\nu_{PD}}F_{SI}}\right)$$
$$= \sigma\sqrt{\frac{U_{obs}}{n-1}}\sqrt{\frac{\Gamma\left(\frac{1}{2}(\nu_{BD}-1)\right)\Gamma\left(\frac{1}{2}(n_{mis}+\nu_{PD})\right)}{\Gamma\left(\frac{\nu_{BD}}{2}\right)\Gamma\left(\frac{1}{2}(n_{mis}+\nu_{PD}-1)\right)}}$$

(36)

And under ML imputation, the distributions are

$$\hat{\mu}_{\infty I,M} = \mu + \sigma\left(\bar{Z}_{obs} + \frac{n_{mis}}{n}\sqrt{\frac{U_{obs}}{n_{obs}+c_M}}E(\bar{Z}_{imp,m})\right)$$
$$= \mu + \sigma\bar{Z}_{obs}$$

$$\hat{\sigma}^2_{\infty I,M} = \sigma^2\frac{U_{obs}}{n-1}\left(1 + \frac{1}{\nu_{obs}+c_M}\left(E(U_{imp,m}) + \frac{n_{obs}}{n}E(n_{mis}\bar{Z}^2_{imp,m})\right)\right)$$
$$= \sigma^2\frac{U_{obs}}{n-1}\left(1 + \frac{1}{n_{obs}+c_M}\left(n_{mis}-1+\frac{n_{obs}}{n}\right)\right)$$

(37)

$$\hat{\sigma}_{\infty I,M} = \sigma \sqrt{\frac{U_{obs}}{n-1}} E\left(\sqrt{1 + \frac{1}{\nu_{obs} + c_M}\left(U_{imp} + \frac{n_{obs}}{n}n_{mis}\bar{Z}_{imp}^2\right)}\right)$$

Note that the distribution of $\hat{\mu}_{\infty I}$ is the same under ML imputation as under PD imputation, but the distributions of $\hat{\sigma}_{\infty I}^2$ and $\hat{\sigma}_{\infty I}$ are different for finite $n_{obs}$. Asymptotically (as $n_{obs} \to \infty$) the distributions of all three $\infty$I estimators—$\hat{\mu}_{\infty I}$, $\hat{\sigma}_{\infty I}^2$, and $\hat{\mu}_{\infty I}$—are the same under ML imputation as under PD imputation.

Note also that the expression for $\hat{\sigma}_{\infty I,M}$ includes an expectation that has no closed-form solution and must be calculated numerically.

# TABLES AND FIGURES

Table 1. Observed-data estimators

a. Bias and standard error

| Parameter | Estimator | Bias | Standard error (SE) |
|:---:|:---:|:---:|:---:|
| $\boldsymbol{\mu}$ | M | $\mathbf{0}$ | $\frac{\sigma}{\sqrt{n_{\text{obs}}}}$ |
| | PD | $0$ | $\sigma\sqrt{\frac{2n_{\text{obs}}+\nu_{\text{prior}}-4}{n_{\text{obs}}(n_{\text{obs}}+\nu_{\text{prior}}-3)}}$ |
| $\boldsymbol{\sigma^2}$ | M | $-\sigma^2\frac{c_M}{c_M+n_{\text{obs}}-1}$ | $\sigma^2\frac{\sqrt{2(n_{\text{obs}}-1)}}{c_M+n_{\text{obs}}-1}$ |
| | PD | $\sigma^2\frac{2-\nu_{\text{prior}}}{n_{\text{obs}}+\nu_{\text{prior}}-3}$ | $\sigma^2\sqrt{\frac{2(n_{\text{obs}}-1)(2n_{\text{obs}}+\nu_{\text{prior}}-4)}{(n_{\text{obs}}+\nu_{\text{prior}}-5)(n_{\text{obs}}+\nu_{\text{prior}}-3)^2}}$ |
| $\boldsymbol{\sigma}$ | M | $\sigma\left(\sqrt{\frac{2}{c_M+n_{\text{obs}}-1}}\frac{\Gamma\left[\frac{n_{\text{obs}}}{2}\right]}{\Gamma\left[\frac{1}{2}(n_{\text{obs}}-1)\right]}-1\right)$ | $\sigma\sqrt{\frac{1}{c_M+n_{\text{obs}}-1}\left(n_{\text{obs}}-1-\frac{2\Gamma\left[\frac{n_{\text{obs}}}{2}\right]^2}{\Gamma\left[\frac{1}{2}(n_{\text{obs}}-1)\right]^2}\right)}$ |
| | PD | $\sigma\left(\frac{\Gamma\left[\frac{n_{\text{obs}}}{2}\right]\Gamma\left[\frac{1}{2}(n_{\text{obs}}+\nu_{\text{prior}}-2)\right]}{\Gamma\left[\frac{1}{2}(n_{\text{obs}}-1)\right]\Gamma\left[\frac{1}{2}(n_{\text{obs}}+\nu_{\text{prior}}-1)\right]}-1\right)$ | $\sigma\sqrt{\frac{n_{\text{obs}}-1}{n_{\text{obs}}+\nu_{\text{prior}}-3}-\frac{\Gamma\left[\frac{n_{\text{obs}}}{2}\right]^2\Gamma\left[\frac{1}{2}(n_{\text{obs}}+\nu_{\text{prior}}-2)\right]^2}{\Gamma\left[\frac{1}{2}(n_{\text{obs}}-1)\right]^2\Gamma\left[\frac{1}{2}(n_{\text{obs}}+\nu_{\text{prior}}-1)\right]^2}}$ |

*b.* Illustrative results

| Number of observations ($n_{obs}$) | Estimator | Common name | $\mu = 1$ | | | | $\sigma^2 = 1$ | | | | $\sigma = 1$ | | | |
|---|---|---|---|---|---|---|---|---|---|---|---|---|---|---|
| | | | Expectation | Bias | SE | RMSE | Expectation | Bias | SE | RMSE | Expectation | Bias | SE | RMSE |
| 5 | M0 | MVU estimator for $\sigma^2$ | 1.00 | 0.00 | 0.45 | 0.45 | 1.00 | 0.00 | 0.71 | 0.71 | 0.94 | -0.06 | 0.34 | 0.35 |
| | M1 | ML estimator | 1.00 | 0.00 | 0.45 | 0.45 | 0.80 | -0.20 | 0.57 | 0.60 | 0.84 | -0.16 | 0.31 | 0.34 |
| | M2 | MMSE estimator for $\sigma^2$ | 1.00 | 0.00 | 0.45 | 0.45 | 0.67 | -0.33 | 0.47 | 0.58 | 0.77 | -0.23 | 0.28 | 0.36 |
| | PD-2 | uniform prior for $\sigma^2$ | 1.00 | 0.00 | undef | undef | undef | undef | undef | undef | 2.36 | 1.36 | undef | undef |
| | PD0 | Jeffreys prior | 1.00 | 0.00 | 0.77 | 0.77 | 2.00 | 1.00 | undef | undef | 1.18 | 0.18 | 0.78 | 0.80 |
| | PD2 | SOUP prior for $\sigma^2$ | 1.00 | 0.00 | 0.63 | 0.63 | 1.00 | 0.00 | 1.41 | 1.41 | 0.88 | -0.12 | 0.47 | 0.48 |
| | PD4 | | 1.00 | 0.00 | 0.58 | 0.58 | 0.67 | -0.33 | 0.75 | 0.82 | 0.74 | -0.26 | 0.35 | 0.44 |
| | PD6 | | 1.00 | 0.00 | 0.55 | 0.55 | 0.50 | -0.50 | 0.50 | 0.71 | 0.64 | -0.36 | 0.29 | 0.46 |
| | PD7 | | 1.00 | 0.00 | 0.54 | 0.54 | 0.44 | -0.56 | 0.43 | 0.70 | 0.61 | -0.39 | 0.27 | 0.47 |
| 20 | M0 | MVU estimator for $\sigma^2$ | 1.00 | 0.00 | 0.22 | 0.22 | 1.00 | 0.00 | 0.32 | 0.32 | 0.99 | -0.01 | 0.16 | 0.16 |
| | M1 | ML estimator | 1.00 | 0.00 | 0.22 | 0.22 | 0.95 | -0.05 | 0.31 | 0.31 | 0.96 | -0.04 | 0.16 | 0.16 |
| | M2 | MMSE estimator for $\sigma^2$ | 1.00 | 0.00 | 0.22 | 0.22 | 0.90 | -0.10 | 0.29 | 0.31 | 0.94 | -0.06 | 0.15 | 0.17 |
| | PD-2 | uniform prior for $\sigma^2$ | 1.00 | 0.00 | 0.34 | 0.34 | 1.27 | 0.27 | 0.66 | 0.72 | 1.09 | 0.09 | 0.27 | 0.29 |
| | PD0 | Jeffreys prior | 1.00 | 0.00 | 0.33 | 0.33 | 1.12 | 0.12 | 0.56 | 0.57 | 1.03 | 0.03 | 0.25 | 0.25 |
| | PD2 | SOUP prior for $\sigma^2$ | 1.00 | 0.00 | 0.32 | 0.32 | 1.00 | 0.00 | 0.49 | 0.49 | 0.97 | -0.03 | 0.23 | 0.23 |
| | PD4 | | 1.00 | 0.00 | 0.31 | 0.31 | 0.90 | -0.10 | 0.43 | 0.44 | 0.93 | -0.07 | 0.21 | 0.22 |
| | PD6 | | 1.00 | 0.00 | 0.30 | 0.30 | 0.83 | -0.17 | 0.38 | 0.42 | 0.89 | -0.11 | 0.20 | 0.23 |
| | PD7 | | 1.00 | 0.00 | 0.30 | 0.30 | 0.79 | -0.21 | 0.36 | 0.42 | 0.87 | -0.13 | 0.19 | 0.23 |
| 100 | M0 | MVU estimator for $\sigma^2$ | 1.00 | 0.00 | 0.10 | 0.10 | 1.00 | 0.00 | 0.14 | 0.14 | 1.00 | 0.00 | 0.07 | 0.07 |
| | M1 | ML estimator | 1.00 | 0.00 | 0.10 | 0.10 | 0.99 | -0.01 | 0.14 | 0.14 | 0.99 | -0.01 | 0.07 | 0.07 |
| | M2 | MMSE estimator for $\sigma^2$ | 1.00 | 0.00 | 0.10 | 0.10 | 0.98 | -0.02 | 0.14 | 0.14 | 0.99 | -0.01 | 0.07 | 0.07 |
| | PD-2 | uniform prior for $\sigma^2$ | 1.00 | 0.00 | 0.14 | 0.14 | 1.04 | 0.04 | 0.21 | 0.22 | 1.02 | 0.02 | 0.10 | 0.10 |
| | PD0 | Jeffreys prior | 1.00 | 0.00 | 0.14 | 0.14 | 1.02 | 0.02 | 0.21 | 0.21 | 1.01 | 0.01 | 0.10 | 0.10 |
| | PD2 | SOUP prior for $\sigma^2$ | 1.00 | 0.00 | 0.14 | 0.14 | 1.00 | 0.00 | 0.20 | 0.20 | 0.99 | -0.01 | 0.10 | 0.10 |
| | PD4 | | 1.00 | 0.00 | 0.14 | 0.14 | 0.98 | -0.02 | 0.20 | 0.20 | 0.99 | -0.01 | 0.10 | 0.10 |
| | PD6 | | 1.00 | 0.00 | 0.14 | 0.14 | 0.96 | -0.04 | 0.19 | 0.20 | 0.98 | -0.02 | 0.10 | 0.10 |
| | PD7 | | 1.00 | 0.00 | 0.14 | 0.14 | 0.95 | -0.05 | 0.19 | 0.20 | 0.97 | -0.03 | 0.10 | 0.10 |

*Note.* "Undef" means that the quantity is undefined, because the formula requires either dividing by zero or taking the square root of a negative number.

Table 2. Single imputation (SI) estimators.

a. Bias and standard error

| Parameter | Estimator | Bias | Standard error (SE) |
|---|---|---|---|
| $\boldsymbol{\mu}$ | M | $\mathbf{0}$ | $\sigma\sqrt{\frac{1}{n_{\text{obs}}}+\frac{n_{\text{mis}}(n_{\text{obs}}-1)}{n^2(n_{\text{obs}}+c_M-1)}}$ |
| | PD | $0$ | $\sigma\sqrt{\frac{1}{n_{\text{obs}}}\left(1+\frac{n_{\text{mis}}(n_{\text{obs}}-1)}{n(n_{\text{obs}}+\nu_{\text{prior}}-3)}\right)}$ |
| $\boldsymbol{\sigma^2}$ | M | $-\sigma^2\dfrac{n_{\text{mis}}(c_M n+n_{\text{obs}}-1)}{(n-1)n(c_M+n_{\text{obs}}-1)}$ | $\sigma^2\frac{\sqrt{2(n_{\text{obs}}-1)}}{(n-1)n(c_M+n_{\text{obs}})}\sqrt{\begin{array}{l}n_{\text{mis}}^4+(2c_M+5n_{\text{obs}}-3)n_{\text{mis}}^3\\+(c_M^2+(6n_{\text{obs}}-4)c_M+n_{\text{obs}}(8n_{\text{obs}}-9)+3\\+n_{\text{obs}}(2c_M^2+6(n_{\text{obs}}-1)c_M+n_{\text{obs}}(5n_{\text{obs}}-9))\cdot\\+n_{\text{obs}}^2(c_M+n_{\text{obs}}-1)^2\end{array}}$ |
| | PD | $\dfrac{\sigma^2}{n-1}\left(\dfrac{n_{\text{mis}}(\nu_{\text{prior}}-2)}{n_{\text{obs}}+\nu_{\text{prior}}-3}\right)$ | $\frac{\sigma^2}{n-1}\sqrt{\frac{2(n_{\text{obs}}-1)}{n_{\text{obs}}+\nu_{\text{prior}}-3}}\sqrt{n+\nu_{\text{prior}}-3+\frac{n_{\text{mis}}(2n_{\text{obs}}+\nu_{\text{prior}}-4)}{(n_{\text{obs}}+\nu_{\text{prior}}-5)(n_{\text{obs}}+\nu_{\text{prior}}-3)}}$ |
| $\boldsymbol{\sigma}$ | M | $\sigma\left(\sqrt{2}\,\frac{\Gamma\left[\frac{n_{\text{obs}}}{2}\right]}{\Gamma\left[\frac{1}{2}(n_{\text{obs}}-1)\right]}\frac{E\left(\sqrt{n_{\text{mis}}n_{\text{obs}}\bar{Z}_{imp}^2+n(c_M+n_{\text{obs}}-1+U_{\text{imp}})}\right)}{\sqrt{(n^2-n)(c_M+n_{\text{obs}}-1)}}-1\right)$ | $\sqrt{E\left(\hat{\sigma}_{SI,M}^2\right)-\left(E\left(\hat{\sigma}_{SI,M}\right)\right)^2}$ |
| | PD | $\sigma\sqrt{\frac{2}{n-1}}\left(\frac{\Gamma\left[\frac{n_{\text{obs}}}{2}\right]\Gamma\left[\frac{1}{2}(n_{\text{obs}}+\nu_{\text{prior}}-2)\right]\Gamma\left[\frac{1}{2}(n+\nu_{\text{prior}}-1)\right]}{\Gamma\left[\frac{1}{2}(n_{\text{obs}}-1)\right]\Gamma\left[\frac{1}{2}(n_{\text{obs}}+\nu_{\text{prior}}-1)\right]\Gamma\left[\frac{1}{2}(n+\nu_{\text{prior}}-2)\right]}-1\right)$ | $\frac{\sigma}{\sqrt{n-1}}\sqrt{\begin{array}{l}\frac{(n_{\text{obs}}-1)(n+\nu_{\text{prior}}-3)}{n_{\text{obs}}+\nu_{\text{prior}}-3}\\-\frac{2\Gamma\left[\frac{n_{\text{obs}}}{2}\right]^2\Gamma\left[\frac{1}{2}(n+\nu_{\text{prior}}-1)\right]^2\Gamma\left[\frac{1}{2}(n_{\text{obs}}+\nu_{\text{prior}}-2)\right]^2}{\Gamma\left[\frac{1}{2}(n_{\text{obs}}-1)\right]^2\Gamma\left[\frac{1}{2}(n+\nu_{\text{prior}}-2)\right]^2\Gamma\left[\frac{1}{2}(n_{\text{obs}}+\nu_{\text{prior}}-1)\right]^2}\end{array}}$ |

b.   Illustrative results with $n_{mis} = n_{obs}$.

| Number of observations s($n_{obs}$) | Estimator | Parameter = true value | | | | | | | | | | | |
| --- | --- | --- | --- | --- | --- | --- | --- | --- | --- | --- | --- | --- | --- |
| | | $\mu = 1$ | | | | $\sigma^2 = 1$ | | | | $\sigma = 1$ | | | |
| | | Expectation | Bias | SE | RMSE | Expectation | Bias | SE | RMSE | Expectation | Bias | SE | RMSE |
| 5 | M0 | 1.00 | 0.00 | 0.50 | 0.50 | 0.94 | -0.06 | 0.78 | 0.78 | 0.90 | -0.10 | 0.36 | 0.38 |
| | M1 | 1.00 | 0.00 | 0.50 | 0.50 | 0.84 | -0.16 | 0.68 | 0.69 | 0.85 | -0.15 | 0.34 | 0.37 |
| | M2 | 1.00 | 0.00 | 0.50 | 0.50 | 0.78 | -0.22 | 0.61 | 0.65 | 0.82 | -0.18 | 0.32 | 0.37 |
| | PD-2 | 1.00 | 0.00 | undef | undef | undef | undef | undef | undef | 1.85 | 0.85 | undef | undef |
| | PD0 | 1.00 | 0.00 | 0.63 | 0.63 | 1.56 | 0.56 | undef | undef | 1.08 | 0.08 | 0.63 | 0.63 |
| | PD2 | 1.00 | 0.00 | 0.55 | 0.55 | 1.00 | 0.00 | 1.15 | 1.15 | 0.91 | -0.09 | 0.42 | 0.43 |
| | PD4 | 1.00 | 0.00 | 0.52 | 0.52 | 0.81 | -0.19 | 0.75 | 0.77 | 0.83 | -0.17 | 0.35 | 0.39 |
| | PD6 | 1.00 | 0.00 | 0.50 | 0.50 | 0.72 | -0.28 | 0.60 | 0.66 | 0.79 | -0.21 | 0.31 | 0.38 |
| | PD7 | 1.00 | 0.00 | 0.49 | 0.49 | 0.69 | -0.31 | 0.56 | 0.64 | 0.77 | -0.23 | 0.30 | 0.38 |
| 20 | M0 | 1.00 | 0.00 | 0.25 | 0.25 | 0.99 | -0.01 | 0.36 | 0.36 | 0.98 | -0.02 | 0.18 | 0.18 |
| | M1 | 1.00 | 0.00 | 0.25 | 0.25 | 0.96 | -0.04 | 0.35 | 0.35 | 0.97 | -0.03 | 0.17 | 0.18 |
| | M2 | 1.00 | 0.00 | 0.25 | 0.25 | 0.94 | -0.06 | 0.34 | 0.35 | 0.95 | -0.05 | 0.17 | 0.18 |
| | PD-2 | 1.00 | 0.00 | 0.29 | 0.29 | 1.14 | 0.14 | 0.51 | 0.53 | 1.04 | 0.04 | 0.22 | 0.23 |
| | PD0 | 1.00 | 0.00 | 0.28 | 0.28 | 1.06 | 0.06 | 0.46 | 0.46 | 1.01 | 0.01 | 0.21 | 0.21 |
| | PD2 | 1.00 | 0.00 | 0.27 | 0.27 | 1.00 | 0.00 | 0.41 | 0.41 | 0.98 | -0.02 | 0.20 | 0.20 |
| | PD4 | 1.00 | 0.00 | 0.27 | 0.27 | 0.95 | -0.05 | 0.38 | 0.39 | 0.96 | -0.04 | 0.19 | 0.19 |
| | PD6 | 1.00 | 0.00 | 0.27 | 0.27 | 0.91 | -0.09 | 0.36 | 0.37 | 0.94 | -0.06 | 0.18 | 0.19 |
| | PD7 | 1.00 | 0.00 | 0.26 | 0.26 | 0.89 | -0.11 | 0.35 | 0.36 | 0.93 | -0.07 | 0.18 | 0.19 |
| 100 | M0 | 1.00 | 0.00 | 0.11 | 0.11 | 1.00 | 0.00 | 0.16 | 0.16 | 1.00 | 0.00 | 0.08 | 0.08 |
| | M1 | 1.00 | 0.00 | 0.11 | 0.11 | 0.99 | -0.01 | 0.16 | 0.16 | 0.99 | -0.01 | 0.08 | 0.08 |
| | M2 | 1.00 | 0.00 | 0.11 | 0.11 | 0.99 | -0.01 | 0.16 | 0.16 | 0.99 | -0.01 | 0.08 | 0.08 |
| | PD-2 | 1.00 | 0.00 | 0.12 | 0.12 | 1.02 | 0.02 | 0.18 | 0.18 | 1.01 | 0.01 | 0.09 | 0.09 |
| | PD0 | 1.00 | 0.00 | 0.12 | 0.12 | 1.01 | 0.01 | 0.18 | 0.18 | 1.00 | 0.00 | 0.09 | 0.09 |
| | PD2 | 1.00 | 0.00 | 0.12 | 0.12 | 1.00 | 0.00 | 0.18 | 0.18 | 1.00 | 0.00 | 0.09 | 0.09 |
| | PD4 | 1.00 | 0.00 | 0.12 | 0.12 | 0.99 | -0.01 | 0.17 | 0.17 | 0.99 | -0.01 | 0.09 | 0.09 |
| | PD6 | 1.00 | 0.00 | 0.12 | 0.12 | 0.98 | -0.02 | 0.17 | 0.17 | 0.99 | -0.01 | 0.09 | 0.09 |
| | PD7 | 1.00 | 0.00 | 0.12 | 0.12 | 0.98 | -0.02 | 0.17 | 0.17 | 0.98 | -0.02 | 0.08 | 0.09 |

*Note.* The expectation of $\widehat{\sigma}_{SI,M}$ was calculated numerically, using the NExpectation function in Mathematica 8.

Table 3. Infinite imputation ($\infty$I) estimators and multiple imputation (MI) estimators

a. Standard errors. (Biases are the same as under single imputation (SI).)

| Parameter | Estimator | Standard error | |
|---|---|---|---|
| | | Infinite imputation ($\infty$I) | Multiple imputation (MI) |
| $\boldsymbol{\mu}$ | M | $\sigma\sqrt{\dfrac{1}{n_{\text{obs}}}}$ | $\sigma\sqrt{\dfrac{1}{n_{\text{obs}}}+\dfrac{n_{\text{mis}}(n_{\text{obs}}-1)}{Dn^2(c_M+n_{\text{obs}}-1)}}$ |
| | PD | $\sigma\sqrt{\dfrac{1}{n_{\text{obs}}}}$ | $\sigma\sqrt{\dfrac{1}{n_{\text{obs}}}\left(1+\dfrac{n_{\text{mis}}(n_{\text{obs}}-1)}{Dn(n_{\text{obs}}+\nu_{\text{prior}}-3)}\right)}$ |
| $\boldsymbol{\sigma^2}$ | M | $\sigma^2\dfrac{\sqrt{2(n_{\text{obs}}-1)}}{n-1}\dfrac{n(c_M+n-2)+n_{\text{obs}}}{n(c_M+n_{\text{obs}}-1)}$ | $\sqrt{\left(1-\dfrac{1}{D}\right)V(\hat\sigma^2_{\infty I,M})+\dfrac{1}{D}V(\hat\sigma^2_{SI,M})}$ |
| | PD | $\sigma^2\dfrac{\sqrt{2(n_{\text{obs}}-1)}}{n-1}\dfrac{n+\nu_{\text{prior}}-3}{n_{\text{obs}}+\nu_{\text{prior}}-3}$ | $\sqrt{\left(1-\dfrac{1}{D}\right)V(\hat\sigma^2_{\infty I,BD})+\dfrac{1}{D}V(\hat\sigma^2_{SI,BD})}$ |
| $\boldsymbol{\sigma}$ | M | $\sigma\,E\left(\sqrt{n(U_{imp}+c_M+n_{\text{obs}}-1)+n_{\text{mis}}n_{\text{obs}}\overline{Z}^2_{imp}}\right)$ $\times\sqrt{\dfrac{n_{\text{obs}}\Gamma\left[\frac{1}{2}(n_{\text{obs}}-1)\right]^2-\Gamma\left[\frac{1}{2}(n_{\text{obs}}-1)\right]^2-2\Gamma\left[\frac{n_{\text{obs}}}{2}\right]^2}{(n-1)n(c_M+n_{\text{obs}}-1)\Gamma\left[\frac{1}{2}(n_{\text{obs}}-1)\right]^2}}$ | $\sqrt{\left(1-\dfrac{1}{D}\right)V(\hat\sigma_{\infty I,M})+\dfrac{1}{D}V(\hat\sigma_{SI,M})}$ |
| | PD | $\dfrac{\sigma}{\sqrt{n-1}}\dfrac{(n+\nu_{\text{prior}}-3)}{(n_{\text{obs}}+\nu_{\text{prior}}-3)}\left(n_{\text{obs}}-1-\dfrac{2\Gamma\left[\frac{n_{\text{obs}}}{2}\right]^2}{\Gamma\left[\frac{1}{2}(n_{\text{obs}}-1)\right]^2}\right)$ | $\sqrt{\left(1-\dfrac{1}{D}\right)V(\hat\sigma_{\infty I,BD})+\dfrac{1}{D}V(\hat\sigma_{SI,BD})}$ |

b. <u>Illustrative results with $n_{mis} = n_{obs}$ and $D$=5 imputations.</u>

| Estimator | Common name | $\mu = 1$ | | | | $\sigma^2 = 1$ | | | | $\sigma = 1$ | | | |
|---|---|---|---|---|---|---|---|---|---|---|---|---|---|
| | | Expectation | Bias | SE | RMSE | Expectation | Bias | SE | RMSE | Expectation | Bias | SE | RMSE |
| M0 | MVU estimator for $\sigma^2$ | 1.00 | 0.00 | 0.46 | 0.46 | 0.94 | -0.06 | 0.69 | 0.69 | 0.90 | -0.10 | 0.32 | 0.34 |
| M1 | ML estimator | 1.00 | 0.00 | 0.46 | 0.46 | 0.84 | -0.16 | 0.61 | 0.63 | 0.85 | -0.15 | 0.31 | 0.34 |
| M2 | MMSE estimator for $\sigma^2$ | 1.00 | 0.00 | 0.46 | 0.46 | 0.78 | -0.22 | 0.56 | 0.60 | 0.82 | -0.18 | 0.30 | 0.35 |
| PD-2 | uniform prior for $\sigma^2$ | 1.00 | 0.00 | undef | undef | undef | undef | undef | undef | 1.85 | 0.85 | undef | undef |
| PD0 | Jeffreys prior | 1.00 | 0.00 | 0.49 | 0.49 | 1.56 | 0.56 | undef | undef | 1.08 | 0.08 | 0.47 | 0.48 |
| PD2 | SOUP prior for $\sigma^2$ | 1.00 | 0.00 | 0.47 | 0.47 | 1.00 | 0.00 | 0.82 | 0.82 | 0.91 | -0.09 | 0.36 | 0.37 |
| PD4 | | 1.00 | 0.00 | 0.46 | 0.46 | 0.81 | -0.19 | 0.61 | 0.64 | 0.83 | -0.17 | 0.32 | 0.36 |
| PD6 | | 1.00 | 0.00 | 0.46 | 0.46 | 0.72 | -0.28 | 0.53 | 0.60 | 0.79 | -0.21 | 0.30 | 0.36 |
| PD7 | | 1.00 | 0.00 | 0.46 | 0.46 | 0.69 | -0.31 | 0.50 | 0.59 | 0.77 | -0.23 | 0.29 | 0.37 |
| M0 | MVU estimator for $\sigma^2$ | 1.00 | 0.00 | 0.23 | 0.23 | 0.99 | -0.01 | 0.33 | 0.33 | 0.98 | -0.02 | 0.16 | 0.16 |
| M1 | ML estimator | 1.00 | 0.00 | 0.23 | 0.23 | 0.96 | -0.04 | 0.32 | 0.32 | 0.97 | -0.03 | 0.16 | 0.16 |
| M2 | MMSE estimator for $\sigma^2$ | 1.00 | 0.00 | 0.23 | 0.23 | 0.94 | -0.06 | 0.31 | 0.32 | 0.95 | -0.05 | 0.16 | 0.16 |
| PD-2 | uniform prior for $\sigma^2$ | 1.00 | 0.00 | 0.24 | 0.24 | 1.14 | 0.14 | 0.40 | 0.42 | 1.04 | 0.04 | 0.18 | 0.19 |
| PD0 | Jeffreys prior | 1.00 | 0.00 | 0.24 | 0.24 | 1.06 | 0.06 | 0.37 | 0.37 | 1.01 | 0.01 | 0.18 | 0.18 |
| PD2 | SOUP prior for $\sigma^2$ | 1.00 | 0.00 | 0.23 | 0.23 | 1.00 | 0.00 | 0.34 | 0.34 | 0.98 | -0.02 | 0.17 | 0.17 |
| PD4 | | 1.00 | 0.00 | 0.23 | 0.23 | 0.95 | -0.05 | 0.32 | 0.33 | 0.96 | -0.04 | 0.16 | 0.17 |
| PD6 | | 1.00 | 0.00 | 0.23 | 0.23 | 0.91 | -0.09 | 0.31 | 0.32 | 0.94 | -0.06 | 0.16 | 0.17 |
| PD7 | | 1.00 | 0.00 | 0.23 | 0.23 | 0.89 | -0.11 | 0.30 | 0.32 | 0.93 | -0.07 | 0.16 | 0.17 |
| M0 | MVU estimator for $\sigma^2$ | 1.00 | 0.00 | 0.10 | 0.10 | 1.00 | 0.00 | 0.15 | 0.15 | 1.00 | 0.00 | 0.07 | 0.07 |
| M1 | ML estimator | 1.00 | 0.00 | 0.10 | 0.10 | 0.99 | -0.01 | 0.14 | 0.14 | 0.99 | -0.01 | 0.07 | 0.07 |
| M2 | MMSE estimator for $\sigma^2$ | 1.00 | 0.00 | 0.10 | 0.10 | 0.99 | -0.01 | 0.14 | 0.14 | 0.99 | -0.01 | 0.07 | 0.07 |
| PD-2 | uniform prior for $\sigma^2$ | 1.00 | 0.00 | 0.11 | 0.11 | 1.02 | 0.02 | 0.15 | 0.15 | 1.01 | 0.01 | 0.08 | 0.08 |
| PD0 | Jeffreys prior | 1.00 | 0.00 | 0.10 | 0.10 | 1.01 | 0.01 | 0.15 | 0.15 | 1.00 | 0.00 | 0.07 | 0.07 |
| PD2 | SOUP prior for $\sigma^2$ | 1.00 | 0.00 | 0.10 | 0.10 | 1.00 | 0.00 | 0.15 | 0.15 | 1.00 | 0.00 | 0.07 | 0.07 |
| PD4 | | 1.00 | 0.00 | 0.10 | 0.10 | 0.99 | -0.01 | 0.15 | 0.15 | 0.99 | -0.01 | 0.07 | 0.07 |
| PD6 | | 1.00 | 0.00 | 0.10 | 0.10 | 0.98 | -0.02 | 0.15 | 0.15 | 0.99 | -0.01 | 0.07 | 0.07 |
| PD7 | | 1.00 | 0.00 | 0.10 | 0.10 | 0.98 | -0.02 | 0.15 | 0.15 | 0.98 | -0.02 | 0.07 | 0.07 |

*Note.* The standard error of $\hat{\sigma}_{\infty I,M}$ is calculated numerically.

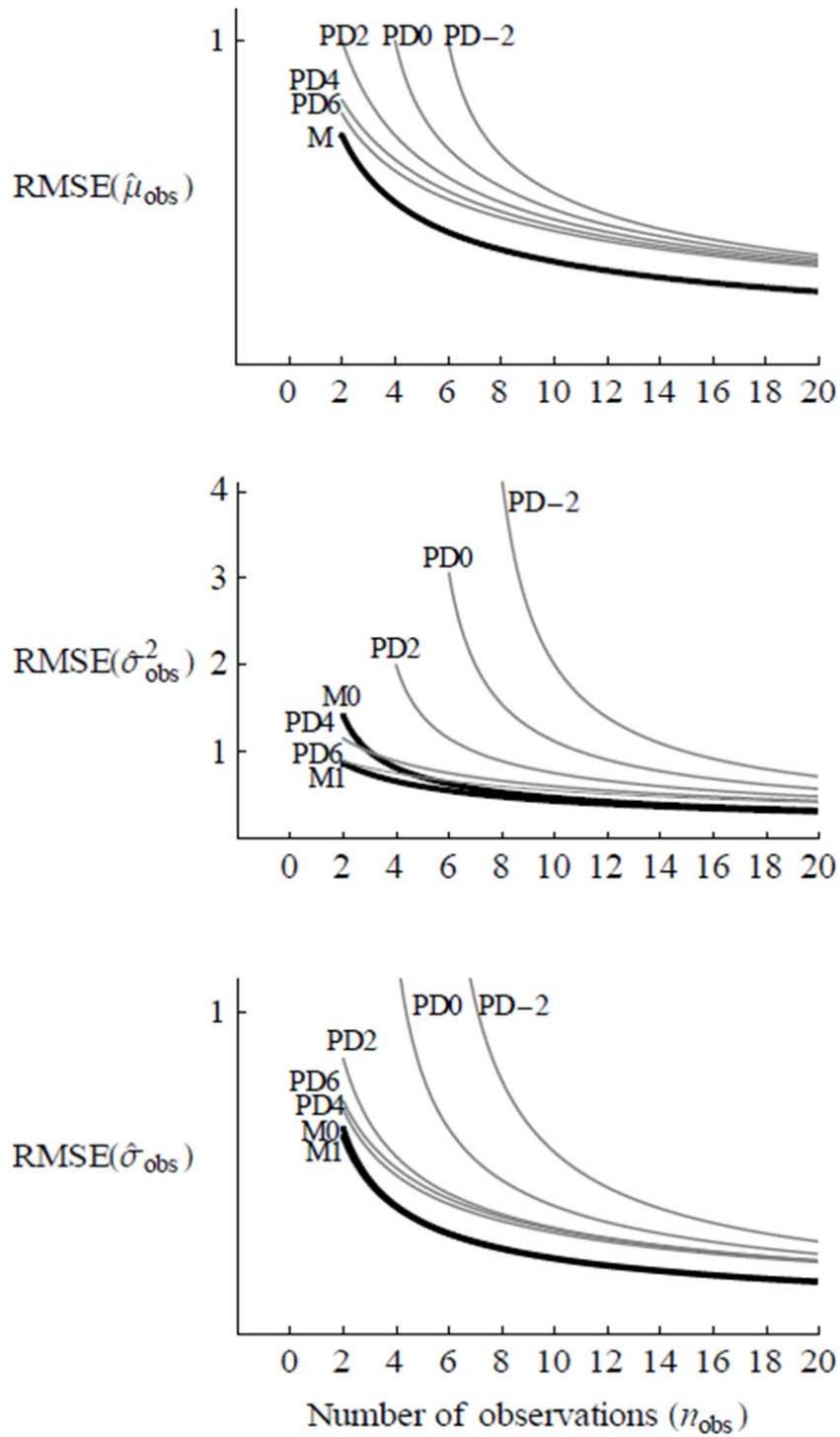

*Figure 1.* Root mean square error (RMSE) for observed-data estimators of the mean, variance, and standard deviation of a normal variable. To standardize the results, the value of $\sigma$ has been set to 1.

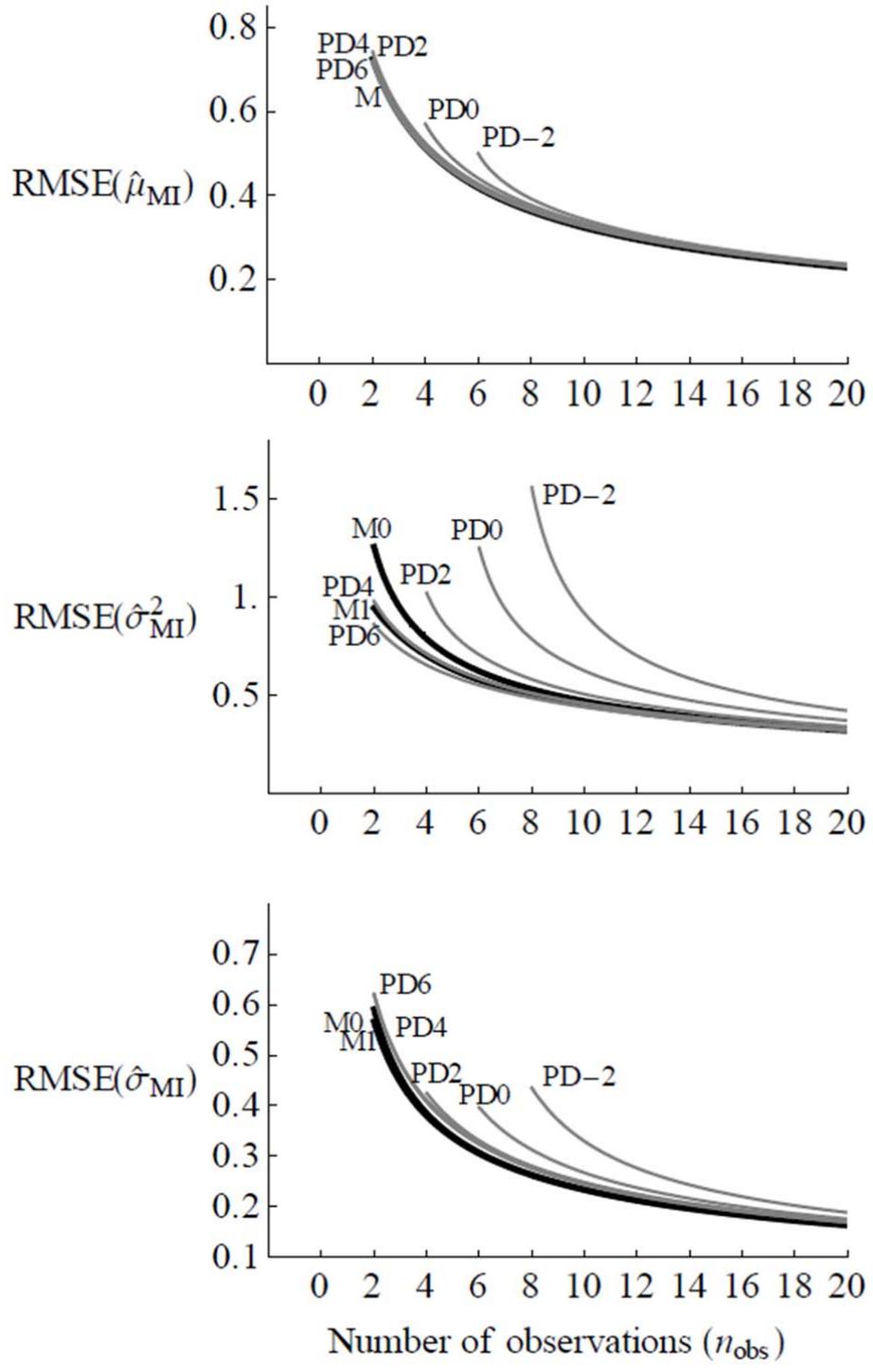

*Figure 2.* RMSE for MI estimators of the mean, variance, and standard deviation of a normal variable. The number of observations $n_{obs}$ increases along the horizontal axis, while the number of imputations is held constant at D=5. To standardize the results, the value of $\sigma$ has been set to 1.

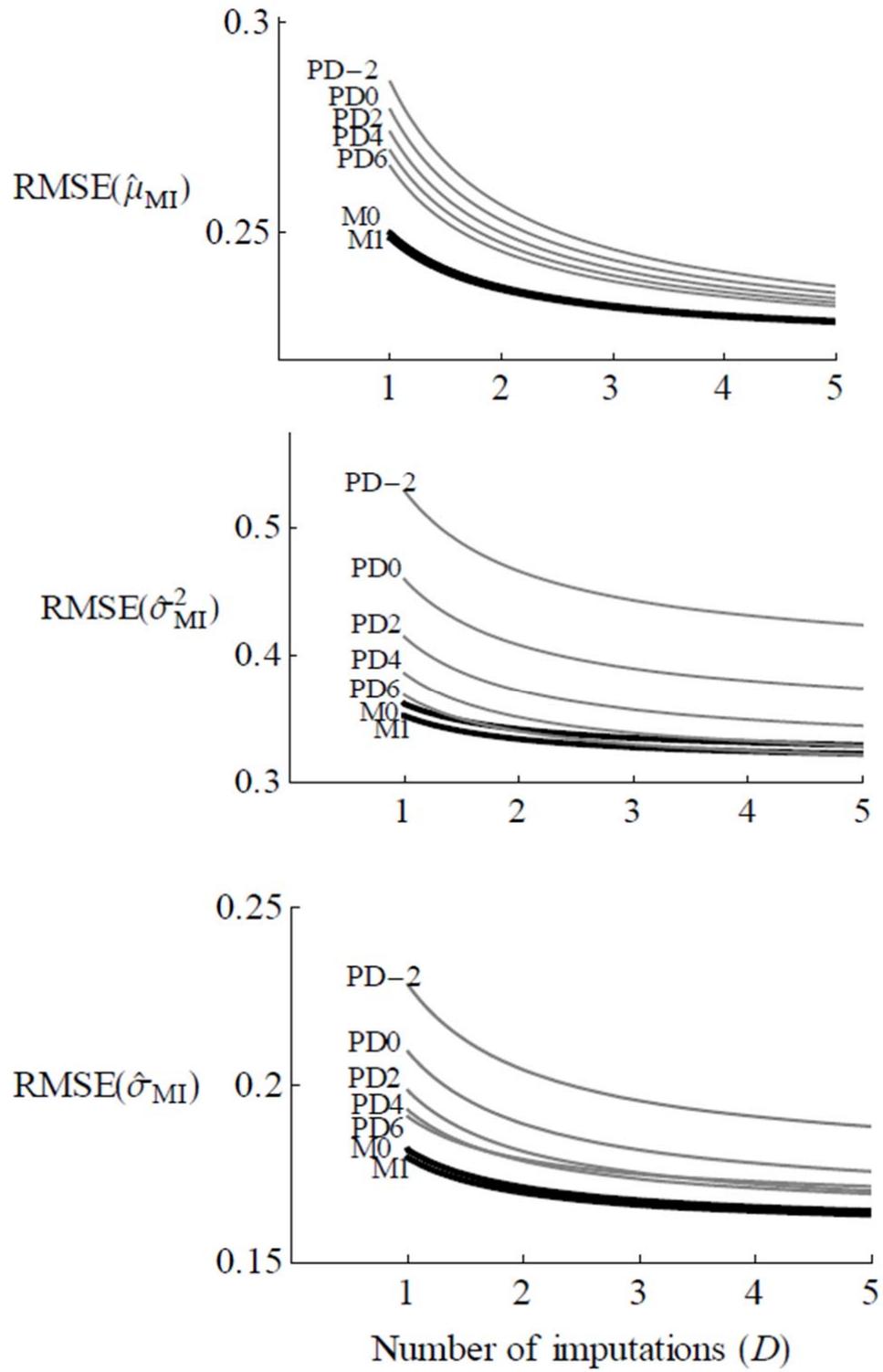

*Figure 3.* RMSE for MI estimators of the mean, variance, and standard deviation of a normal variable. The number of imputations $D$ increases along the horizontal axis, while the number of observed and missing values is held constant at $n_{obs} = n_{mis} = 20$. To standardize the results, the value of $\sigma$ has been set to 1.